\documentclass[12pt,a4paper]{amsart}
\usepackage{amssymb,amsfonts,amsthm,amsmath}
\usepackage{graphicx} 

\theoremstyle{plain}
\newtheorem{theorem}{Theorem}

\numberwithin{equation}{section} \numberwithin{theorem}{section}
\numberwithin{lemma}{section} \numberwithin{definition}{section}
\numberwithin{corollary}{section} \textheight =24cm
\begin{document}
\title{A Diophantine sum with factorials}
\author{Geoffrey B Campbell}
\address{Mathematical Sciences Institute \\
         The Australian National University \\
         ACT, 0200, Australia}
 \email{Geoffrey.Campbell@anu.edu.au}

\author{Aleksander Zujev}
\address{Physics Department \\
         University of California \\
         Davis, California 95616, USA}
\email{azujev@ucdavis.edu}

\begin{abstract}
We give solutions of a Diophantine equation containing factorials, which can be written as a cubic form, or as a sum of binomial coefficients. We also give some solutions to higher degree forms and relate some solutions to an unsolvable Fermat Last Theorem equation.
\end{abstract}
\maketitle

\section{Context of the problem} \label{S:context}
The problem  $x^3+y^3+z^3 = x+y+z$ is listed in the book by
Guy (see ~\cite[Section D8, A pyramidal diophantine equation]{rG1981}) as an unsolved problem of circa 1981.
However, later Guy also mentions Bremner, ~\cite{aB1977} remarking that he "has effectively determined all parametric solutions".
Bremner's paper is thorough and sophisticated, whereas, Oppenheim ~\cite{aO1966} also considers this same problem and provides small solutions.
In volume three of classical Number Theory compendium by Dickson ~\cite{aO1966} the equation $\binom{a}{3} + \binom{b}{3} = \binom{c}{3}$
is simply mentioned as having known solutions, but without examples nor any further remarks.
The problem itself of solving the Diophantine equation $x^3+y^3+z^3 = x+y+z$ was used in a proof of Fermat's Last Theorem by van der Poorten ~\cite{vdP1975}
in the case of index 3, and it turns out equivalent to solving $a!/(a-3)! + b!/(b-3)! = c!/(c-3)!$.

Diophantine equations typically take the simplistic form, for integers $a_j$, $a_k$, $x$, $y$, and fixed positive integers $m$ or $s$,
\begin{equation} \label{E:1.1}
              \sum_{j=0}^{n} a_j x^m = \sum_{k=0}^{r} a_k y^s.
\end{equation}

In this note we consider the form for positive integers $n$, $a-n$, $b-n$, $c-n$,
\begin{equation} \label{E:1.2}
                   a!/(a-n)! + b!/(b-n)! = c!/(c-n)!
\end{equation}

and show that for $n=3$ there are an infinite number of integer solutions, and that, since
\begin{equation} \label{E:1.3}
                   (a+1)!/(a-2)!= a^3  - a
\end{equation}
these reduce to $x^3+y^3+z^3 = x+y+z$ which is in the traditional form (\ref{E:1.1}). It is easy to find solutions of (\ref{E:1.2}) such as the following:-
\begin{equation} \label{E:1.4}    5!/2! + 5!/2! = 6!/3!        \end{equation}
\begin{equation} \label{E:1.5}    10!/7! + 16!/13! = 17!/14!   \end{equation}
\begin{equation} \label{E:1.6}    22!/19! + 56!/53! = 57!/54!  \end{equation}
\begin{equation} \label{E:1.7}    32!/29! + 57!/54! = 60!/57!  \end{equation}
\begin{equation} \label{E:1.8}    41!/38! + 72!/69! = 76!/73!  \end{equation}

Hence, using (\ref{E:1.3}) on (\ref{E:1.4}) to (\ref{E:1.8}) we find solutions of $a^3-a+b^3-b=c^3-c$  for positive integers $a$, $b$ and $c$, for instance,

\begin{equation} \label{E:1.9}    4^3-4+4^3-4=5^3-5;        \end{equation}
\begin{equation} \label{E:1.10}   9^3-9+ 15^3-15=16^3-16;   \end{equation}
\begin{equation} \label{E:1.11}   21^3-21+55^3-55=56^3-56;  \end{equation}
\begin{equation} \label{E:1.12}   31^3-31+56^3-56=59^3-59;  \end{equation}
\begin{equation} \label{E:1.13}   40^3-40+71^3-71=75^3-75.  \end{equation}

We know from Fermat's Last Theorem's that $a^3+b^3=c^3$ is not possible for positive integers $a$, $b$ and $c$, so (\ref{E:1.9}) to (\ref{E:1.13})
is a kind of surprising variant on this. Let's check for solutions to (\ref{E:1.2}) where $n=4$. The equation we are dealing with is therefore,
\begin{equation} \label{E:1.14}a!/(a-4)! + b!/(b-4)! = c!/(c-4)! \end{equation}
or equivalently, for positive integers $a$, $b$ and $c$,
\begin{equation} \label{E:1.15}a(a-1)(a-2)(a-3) + b(b-1)(b-2)(b-3) = c(c-1)(c-2)(c-3). \end{equation}
A shallow search of natural numbers $a$, $b$ and $c$, less than 100 yields no non-trivial cases of (\ref{E:1.15}).
However, the equations $a^3+a+b^3+b=c^3+c-6^3$  and $a^3+a+b^3+b=c^3+c-10^3$ have solutions, some examples as follows,
\begin{equation} \label{E:1.16}1^3+1+8^3+8=9^3+9-6^3;         \end{equation}
and also:
\begin{equation} \label{E:1.17}   1^3+1+9^3+9=12^3+12-10^3,   \end{equation}
\begin{equation} \label{E:1.18}   9^3+9+16^3+16=18^3+18-10^3, \end{equation}
\begin{equation} \label{E:1.19}   8^3+8+22^3+22=23^3+23-10^3. \end{equation}

Evidently this factorial form of problem is relatively new in the ancient literature, as in since 2002.
There appears to be a nontrivial theory worth exploring for it. For example, we show below that there are an infinite
number of solutions to the equation $a!/(a-3)! + b!/(b-3)! = c!/(c-3)!$ for positive integers $a$, $b$ and $c$.
This is equivalent to asserting that there are an infinite number of solutions to the equation $a(a+1)(a+2) + b(b+1)(b+2) = c(c+1)(c+2)$.
There was a brief unresolved online discussion on equations of kind (\ref{E:1.2}) where $n=3$ in 2002.
The discussion thread unsuccessfully sought an identity for these and/or a proof that an infinite number of solutions exist.
We next provide that proof and identity.

\section{Proof of infinite number of solutions} \label{S:InfiniteSolutions}

We start with the
\begin{theorem}
  The equation $a^3-a+b^3-b=c^3-c$ has an infinite number of positive integer solutions.
\end{theorem}

PROOF. Let $a=u-v, b=kv, c=u+v$.  After substitution this becomes $(k^3 - 2)v^2 - 6u^2 = k - 2$,
which is a Pell type equation if we fix a value of $k$.

For $k = 1$; we arrive at $v^2 + 6u^2 = 1$, whose only integer solutions are $u=0$, $v=1$, or $v=-1$. For $k=2$, we have next $v^2 = u^2 + 1$,
whose only integer solutions are also $u=0$, $v=1$, or $v=-1$.

For $k = 3$ we have $25v^2-6u^2=1$,
 which is a Pell equation with an infinite number of solutions.
 END OF PROOF.

Incidentally, for $25v^2-6u^2=1$, the general solution is given by

 $u = ±((5-2 \sqrt{6})^n-(5+2 \sqrt{6})^n)/(2 \sqrt{6})$,

 $5v = ±1/2 (-(5-2 \sqrt{6})^n-(5+2 \sqrt{6})^n)$,

 for non-negative integers $n$.
So we have an infinite family of solutions starting with:

$v=1, u=2, a=1, b=3, c=3$ (trivial solution), then

$v=97, u=198, a=101, b=291, c=295$
...

For $k=4$: $31v^2-3u^2 = 1$

This Pell-type equation has solutions:
$v=14, u=45, a=31, b=56, c=59$

$v=340214, u=1093635, a=753421, b=1360856, c=1433849$
...

For $k=5,6$  there are no solutions.

For $k=7$ we arrive at $341v^2-6u^2 = 5$.

This has an infinite number of solutions, starting with:

$v=13, u=98, a=85, b=91, c=111$

$v=941, u=7094, a=6153, b=6587, c=8035$
...

Any one of the above solvable Pell equations proves the theorem.

\section{Solutions of equation $a!/(a-n)! + b!/(b-n)! = c!/(c-n)!$ for $n=6$} \label{S:Solutions_n6}

We state the
\begin{theorem}
  The equation
  \begin{equation} \label{E:3.1}   a!/(a-6)! + b!/(b-6)! = c!/(c-6)!,   \end{equation}
has positive integer solutions.
\end{theorem}

After a fairly extensive search it seems likely that the only solutions of (\ref{E:3.1}) are:
\begin{equation} \label{E:3.2}   11!/5! + 11!/5! = 12!/6!,   \end{equation}
\begin{equation} \label{E:3.3}   14!/8! + 15!/9! = 16!/10!, \end{equation}
\begin{equation} \label{E:3.4}   19!/13! + 19!/13! = 21!/15!. \end{equation}

There seems no literature related to (\ref{E:3.1}), which could be restated in terms of binomial coefficients as

\begin{equation} \label{E:3.5}   \binom{a}{6} + \binom{b}{6} = \binom{c}{6}   \end{equation}.

Hence, equations (\ref{E:3.2}), (\ref{E:3.3}) and (\ref{E:3.4}) above could be written as solutions to (\ref{E:3.5}) as:

\begin{equation} \label{E:3.6}   \binom{11}{6} + \binom{11}{6} = \binom{12}{6},   \end{equation}
\begin{equation} \label{E:3.7}   \binom{14}{6} + \binom{15}{6} = \binom{16}{6}, \end{equation}
\begin{equation} \label{E:3.8}   \binom{19}{6} + \binom{19}{6} = \binom{21}{6}. \end{equation}

It would perhaps be interesting to solve equations for positive integers $n$ generally,

\begin{equation} \label{E:3.9}   {\binom{a}{n}}^2 + {\binom{b}{n}}^2 = {\binom{c}{n}}^2,   \end{equation}

or for the generalization,

\begin{equation} \label{E:3.10}   {\binom{a}{n}}{\binom{b}{n}} + {\binom{c}{n}}{\binom{d}{n}} = {\binom{e}{n}}{\binom{f}{n}}.   \end{equation}

\end{document}